\documentclass[11pt]{amsart}

\usepackage{amssymb}
\usepackage[colorlinks=true, urlcolor=rltblue, citecolor=drkgreen] {hyperref}
\usepackage{color}
\usepackage{pdfsync}
\definecolor{rltblue}{rgb}{0,0,0.4}
\definecolor{drkgreen}{rgb}{0,0.4,0}
\usepackage{pdfsync}

\newtheorem{thm}{Theorem}[section]
\newtheorem{lemma}[thm]{Lemma}

\newtheorem{theorem}[thm]{Theorem}
\newtheorem{corollary}[thm]{Corollary}

\theoremstyle{definition}
\newtheorem{definition}[thm]{Definition}

\theoremstyle{remark}

\newtheorem{question}{Question}

\theoremstyle{plain}

\newcounter{contenumi}

\def\geqt{\geq_T}

\def\O{{\mathcal O}}

\def\upto{\mathop{\upharpoonright}}
\def\la{\langle}
\def\ra{\rangle}
\def\and{\mathrel{\&}}

\def\isom{\cong}
\def\Si{\Sigma}

\newcommand\rightdate[1]{\footnotetext{  Saved: #1 \\ Compiled: \today}}
\def\A{\mathcal{A}}
\def\B{\mathcal{B}}

\def\E{\mathcal{E}}

\def\om{\omega}
\def\bbar{\bar{b}}

\def\si{\sigma}
\def\b{\beta}

\def\a{\alpha}

\def\P{\mathop{\mathcal P}}

\def\M{{\mathcal M}}
\def\Ztwo{\revmathfont{Z}$_2$}

\def\bfx{\mathbf{x}}

\def\A{\mathcal A}

\def\H{\mathcal H}
\def\xbar{\bar{x}}

\def\ybar{\bar{y}}
\def\abar{\bar{a}}

\def\ctt{{\mathtt c}}

\def\Sic{\Si^\ctt}
\def\Pic{\Pi^\ctt}

\def\g{\gamma}

\newcommand{\revmathfont}[1]{{\textsf{#1}}}

\def\Si{\Sigma}
\def\Ztwo{\revmathfont{Z}$_2$}

\def\Ztwo{\revmathfont{Z}$_2$}

\def\ZFC{\revmathfont{ZFC}}

\def\zs{0^{\#}}
\def\Mo{M}
\def\X{{\mathcal X}}

\def\dtp{\mbox{-}tp}

\def\pbar{\bar{p}}
\def\qbar{\bar{q}}
\def\rbar{\bar{r}}

\def\Sp {\mathop{\mathrm{Sp}}}
\def\om{\omega}

\def\Sp{Sp}

\def\PP{\mathbb{P}}
\def\forces{\Vdash}
\def\bfa{{\bf a}}
\def\eqSIN{\equiv_{I}}
\def\leqSIN{\leq_{I}}
\def\Deltac{\Delta^c}
\def\HF{\mathbb{HF}}
\def\On{On}
\def\Force{Force}

\title{A fixed point for the jump operator on structures}

\author{Antonio Montalb\'an}
\thanks{The author was partially supported by NSF grant DMS-0901169 and the Packard Fellowship.}

\address{Department of Mathematics\\
University of Chicago\\
5734 S. University ave.\\
Chicago, IL 60637, USA}

\email{antonio@math.uchicago.edu}
\urladdr{\href{http://www.math.uchicago.edu/~antonio/index.html}{www.math.uchicago.edu/$\sim$antonio}}

\begin{document}

\rightdate{June 5, 2011 -- submitted}
\maketitle


\begin{abstract}
Assuming that $0^\#$ exists, we prove that there is a structure that can effectively interpret its own jump.
In particular, we get a structure $\mathcal A$ such that 
\[
Sp({\mathcal A}) = \{{\bf x}':{\bf x}\in Sp ({\mathcal A})\},
\]
where $Sp ({\mathcal A})$ is the set of Turing degrees which compute a copy of $\mathcal A$.

It turns out that, more interesting than the result itself, is its unexpected complexity.
We prove that higher-order arithmetic, which is the union of full $n$th-order arithmetic for all $n$, cannot prove the existence of such a structure.
\end{abstract}

\section{Introduction}

Informally, the jump of a abstract structure $\A$ is another structure, $\A'$, obtained by adding to $\A$ relations that code all the $\Si_1$ information about $\A$.

\begin{definition}
For a language $L$, let $\{\varphi_i:i\in\om\}$ be a computable enumeration of all the computably infinitary $\Si_1$ $L$-formulas.
Given an $L$-structure $\A$, let $\A'$ be the structure obtained by adding to $\A$ infinitely many relations $K_i$, for $i\in \om$, where $ \A\models K_i(\xbar)\leftrightarrow \varphi_i(\xbar)$, and where the arity of $K_i$ is the same as the one of $\varphi_i$.
(All the languages we consider in this paper are at most countable and have a computable arity function.)

The computably infinitary $L$-formulas were introduced by Ash (see \cite[Chapter 7]{AK00}).
We use $\Sic_n$ to denote the class of computably infinitary $\Si_n$ formulas.
\end{definition}

The notion of the jump of an abstract structure was introduced recently, independently by various authors.
The definition above is the one that appeared in \cite{MonJumpStr, MonND}.
The other definitions are due to: Baleva \cite{Bal06}, using Moschovakis extensions and a complete $\Sic_1$-relation; A. Soskova and Soskov \cite{Sos07, SS09}, using Moschovakis extensions and coding the forcing relation for $\Pi_1$ formulas; and Stukachev \cite{Stu10}, using hereditarily finite extensions in the context of $\Sigma$-reducibility for structures of arbitrary size.
All these definitions are equivalent in a strong sense, namely up to $\Si$-reducibility.
Puzarenko \cite{Puz09} had also independently introduced an equivalent definition of jump, but did it only for admissible sets rather than general structures.
He was extending a previous notion of jump due to Morozov \cite{Mor04}, that works only for recursively listed admissible sets.

The equivalence of all these definitions reaffirms the naturality of the notion.
A more important reason why this notion is interesting is that it helps explain what is behind many old results in effective structure theory.
For instance, Downey and Jockusch's result \cite{DJ94} that every low Boolean algebra has a computable copy is based on a lemma that can be restated as follows: If $\B$ is a Boolean algebra, and  $0'$ computes a copy of $\B'$, then $\B$ has a computable copy.
Another examples is the result by Ash, Knight, Mennasse and Slaman \cite{AKMS89} and Chisholm \cite{Chi90} that says that every relatively intrinsically  $\Si^0_n$ relation on a structure $\A$ is $\Sic_n$-definable.
One can prove this result by proving only the case $n=1$, where the forcing is very simple, and then applying it to the $(n-1)$st jump of $\A$.
(See \cite[Theorem 3.9]{MonND} for a more detailed proof of this last example. See \cite{MonJumpStr} for more examples.)

Once we are convinced this is a natural notion of jump, the question of whether it is actually a ``jump'' immediately pops up:

\begin{question}  \label{Q: main}
Is there a structure that is equivalent to its own jump?
\end{question}

This question was asked in print in \cite[Remark 1, page 3]{Stu10} and in \cite[Section 7, Question 4]{Puz09}. (The reference in \cite[Remark 1, page 3]{Stu10} to a positive solution for $\equiv_w$ turned out to be incorrect.)

For other jump operators, like the original one on the Turing degrees, or the one on the enumeration degrees, the proofs that there are no fixed points are usually done by a simple diagonalization argument.
This is not the case for the jump of abstract structures.
Furthermore, we will see that, even more interesting than the answer of Question \ref {Q: main} itself, is its unexpected complexity.

For Question \ref {Q: main} to be concrete, we need to specify a notion of equivalence between structures.
(For a study of different notions of equivalence between structures see \cite{Stu07, Kal09}.)
Here is our first candidate.

\begin{definition}
Given two structures $\A$ and $\B$, we say that $\A$ is {\em Muchnik reducible} to $\B$, and write $\A\leq_w \B$, if
\[
\forall X\subseteq \om,\quad X\mbox{ computes a copy of } \B \quad\implies\quad X\mbox{ computes a copy of } \A,
\]
or equivalently, if $\Sp(\B)\subseteq\Sp(\A)$, where $\Sp(\A)$, the {\em degree spectrum of $\A$}, is the set of Turing degrees which compute a copy of $\A$.

This reduction defines a pre-ordering on the class of all countable structures, and hence an equivalence, $\equiv_w$, as usual.
\end{definition}

Even though this notion is not always used as a reducibility, it is widely accepted as a way to measure the computability theoretic complexity of a structure.
The spectrum of a structure behaves well with the notion of jump:
It was proved by A. Soskova and  Soskov \cite{SS09}, and by Montalb\'an \cite{MonJumpStr} independently, that, for every structure $\A$,
\[
\Sp(\A') = \{\bfx': \bfx\in \Sp(\A)\}.
\] 

So, Question \ref {Q: main} for Muchnik equivalence reduces to the question of whether there exists a structure $\A$ for which $\Sp(\A)=\{\bfx': \bfx\in \Sp(\A)\}$.
Here is our first main theorem.

\begin{theorem}(\ZFC + ``$\zs$ exists'') \label{cor: main}
There is a structure $\A$ such that
\[
\Sp(\A) = \Sp(\A').
\]
\end{theorem}
We will prove this theorem in Section \ref {se: fixed point}.

It is well known that there are upward-closed classes of Turing degrees $S$ such that $S=S'$, where $S'$ denotes the set $\{\bfx': \bfx\in S\}$.
For instance, if we take a sequence $\{\bfa_i:i\in\om\}$ of Turing degrees such that $\bfa_i\geqt \bfa_{i+1}'$ for all $i$, then $S=\{\bfx: \exists i\ (\bfx\geqt\bfa_i)\}$ satisfies $S=S'$.
Such a sequence of degrees can be found below Keene's $\O$ (Harrison \cite{Har68}) and its existence can be proved using {\em Arithmetic Transfinite Recursion}.
However, this set $S$ we just constructed is not the spectrum of any structure:
Richter \cite{Ric77} had proved that a spectrum can never be a countable union of upper cones in the Turing degrees.

It is also known that any set $S$ of Turing degrees with $S=S'$ has to be somewhat complex: if $\bfx$ belongs to such a set $S$, then $\bfx$ computes all hyperarithmetic sets (Enderton and Putnam \cite{EP70}).
This level of complexity is still very small compared to what we need to build a structure whose spectrum is jump invariant.

Our proof of Theorem \ref{cor: main} uses $\zs$, which cannot be shown to exist in \ZFC.
The structure $\A$ will be an ill-founded $\om$-model of \ZFC+\revmathfont{V=L}.
Notice that for every $\a\in\A$ which is an ordinal of $\A$, a complete $\Sic_1$ relation on $(L_\a)^\A$ is coded in $(L_{\a+1})^\A$, which implies that $\A$ can effectively interpret the jump of the structure $(L_\a)^\A$.
Thus, it would be enough to build $\A$ so that for some $\a$, which is an ordinal of $\A$, we have $(L_\a)^\A\isom \A$.
In Subsection \ref {ss: fixed point}, we will define such a structure $\A$ as the Skolem hull of a sequence of ordinal indiscernibles, $H$, of type $\zs$ (i.e.\ $\A\models\varphi(\a_1,...,\a_k)$ for $\a_1<...<\a_k\in H$ if and only if $\ulcorner\varphi(x_1,...,x_k)\urcorner\in \zs$), and which has order-type $\om\cdot\om^*$ (i.e.\ $H\isom \cdots+\om+\om+\om$).

We do not know whether there is a proof in \ZFC\ of Theorem \ref {cor: main}.
However, we were able to show that a proof of Theorem \ref{cor: main} has to use techniques that are almost never used in classical mathematics except by logicians.
This implies that a structure whose jump has the same spectrum as itself  will, most likely, never occur naturally outside logic. 
More concretely, we will show that Theorem \ref {cor: main} is not provable in higher-order arithmetic.
Higher-order arithmetic is the union, for all $n\in\om$, of $n$th-order arithmetic which includes the full $n$th-order comprehension scheme (see Section \ref {ss: HOA}).
It is intended to describe the structure $(\om, \P(\om),\P(\P(\om)),...; 0,1,+,\times,<, \in)$.
A great fragment of classical mathematics can be stated and proved in second-order arithmetic.
The small fragment of classical mathematics that cannot be stated in second-order arithmetic might need third-, or at most fourth-order arithmetic, but rarely more, unless it comes from set theory or other parts of logic.

We let \Ztwo\ denote the theory of full second-order arithmetic.
Here is our second main theorem.

\begin{theorem}(\Ztwo)   \label{thm: omega arithmetic}
The existence of a structure $\A$ with $\Sp(\A) = \Sp(\A')$ implies the consistency of higher-order arithmetic.
\end{theorem}

\begin{corollary}\label {cor: omega arithmetic}
Higher-order arithmetic cannot prove that there exists a structure $\A$ with $\Sp(\A) = \Sp(\A')$.
\end{corollary}

To prove Theorem \ref {thm: omega arithmetic}, we need to build a model of higher-order arithmetic out of the structure $\A$.
This model is defined as follows: 
For the first order part, $N$, we use the standard model of the natural numbers.
For $\P(N)$ we take the collection of all the sets $X\subseteq\om$ which are c.e.\ in all copies of $\A$.
For $\P^n(N)$ we use the the collections of all hereditarily countable families $F\in\P^{n}(\om)$ which are uniformly computably enumerable in all copies of $\A$, where $F$ is {\em uniformly computably enumerable} in $Y\subseteq \om$ if there exists a $Y$-c.e.\ set $W$ such that 
\[
F=\{\{...\{\{\{i_0\in\om: \la i_0,i_1,....,i_{n-1}\ra\in W\}: i_1\in\om\}: i_2\in\om\}:\dots\}: i_{n-1}\in\om\}.
\]
The first-order part of this model is standard, so all the effort goes into showing that it satisfies comprehension for sets of all types and formulas of all orders.

In \cite{StuEMU}, Question \ref{Q: main} was asked for $\Sigma$-equivalence, a notion of equivalence between structures that is much stronger than Muchnik equivalence. 
The notion of $\Sigma$-equivalence between abstract structures was introduced and studied in detail in \cite{Khi04, Stu07, Stu08}.
Our proof of Theorem \ref{cor: main} does answer this question, as it actually shows, in \ZFC + ``$\zs$ exists,'' that there is a structure $\A$ such that $\A \equiv_\Sigma\A'$.
Furthermore, we will introduce an even stronger notion of equivalence, $\eqSIN$, and show that our structure satisfies $\A \eqSIN\A'$.
Our new definition of equivalence was inspired by the notion of $\Sigma$-reducibility.

\begin{definition}
Let $\A$ be an $L$-structure, and $\B$ be any structure.
Let us assume that $L$ is a relational language $L=\{P_0,P_1,P_2,...\}$ where $P_i$ has arity $a(i)$; so $\A= (A; P_0^\A,P_1^\A,...)$ and $P_i^\A\subseteq A^{a(i)}$.

We say that $\A$ is {\em effectively interpretable} in $\B$, and write $\A\leqSIN\B$, if, for some $n\in\om$, in $\B$, there is a $\Sic_1$-definable set $D\subseteq B^n$, a $\Deltac_1$-definable relation $\eta\subseteq B^n\times B^n$ which is an equivalence relation on $D$, and  a uniformly $\Deltac_1$-definable sequence of sets $R_i\subseteq B^{n\cdot a(i)}$, closed under the equivalence $\eta$ within $D$, such that 
\[
(A; P_0^\A,P_1^\A,...) \isom (D/\eta; R_0,R_1,...).
\]
The sets $R_i$ do not need to be subsets of $D^{a(i)}$, and, when we refer to the structure $(D/\eta; R_0,R_1,...)$ we, of course, mean $(D/\eta; (R_0\cap D^{a(0)})/\eta,(R_1\cap D^{a(1)})/\eta,...)$.
By {\em  uniformly $\Deltac_1$-definable}, we mean that there is a computable sequence $\{\theta_i: i\in \om\}$ of $\Sic_1$ formulas, and a computable sequence $\{\psi_i: i\in\om\}$ of $\Pic_1$ formulas such that, for all $i$, $\B\models R_i(\xbar)\leftrightarrow \theta_i(\xbar) \leftrightarrow \psi(\xbar)$.
\end{definition}

It is not hard to see that effective interpretability implies Muchnik reducibility.
Furthermore, for the readers familiar with $\Sigma$-reducibility, it is not hard to see that effective interpretability implies $\Sigma$-reducibility too.
This is because a relation $R\subseteq B^m$ is $\Sic_1$ definable in $\B$ if and only if it is $\Si$-definable in $\HF(\B)$ (as proved by Va\v{\i}tsenavichyus; see \cite[Theorem 1]{StuEMU}).
These two reducibilities  are not equivalent: for example, if $(\om;\ )$ is an infinite structure on an empty language, then $(\om;Succ)\leq_\Si (\om;\ ) $ but $(\om,Succ)\not\leqSIN (\om;\ ) $.

Instead of proving Theorem \ref {cor: main}, we will directly prove the following theorem.

\begin{theorem}(\ZFC + ``$\zs$ exists'')    \label{thm: main}
There is a structure $\A$ such that $\A \eqSIN \A'$.
\end{theorem}

\section{The fixed point structure}  \label{se: fixed point}

\subsection{Background on $L$-indiscernibles}

We follow Devlin's book \cite[Chapter 5]{Dev84}.
The reader familiar with $\zs$ can skip this subsection and move to \ref {ss: fixed point}.

We say that a set $S$ of formulas in the language of set theory is an {\em Ehrenfeucht-Mostowski set (E-M set)} if there is a model $\M$ of \ZFC+\revmathfont{V=L} and an infinite set of ordinal indiscernibles $H$ in $\M$ such that $S$ is the set of formulas that are true in $\M$ about increasing tuples from $H$
(i.e., $S=\{\varphi(x_1,...,x_k): \M\models \varphi(\a_1,...,\a_k) \mbox{ for some (all) }\a_1<...<\a_k\in H\}$).
Devlin uses models of the weaker system \revmathfont{BSL+V=L}, instead of \ZFC+\revmathfont{V=L},  in his definition of Ehrenfeucht-Mostowski sets.
However, for our purposes we do not need to deal with \revmathfont{BSL}.

Recall that \revmathfont{V=L} implies that there exists a definable well-ordering $<_L$ of the universe.
We can use this ordering to add Skolem terms $t_\varphi$ to the language for all formulas $\varphi$, and add axioms saying that 
``$t_\varphi(x_1,...,x_n)$ is the $<_L$-least witness for $\varphi(x_1,...,x_n, x)$, if such a witness exists.''
This is an inessential extension of \ZFC+ \revmathfont{V=L}.

From the Ehrenfeucht-Mostowski theorem,  we get that given an E-M set $S$, and a linear ordering $\X$, there is a model $\Mo(S,\X)$ with a set of ordinal indiscernibles $H$, unique up to isomorphism, such that $H$ has order type $\X$ and $\Mo(S,\X)$ is the Skolem hull of of $H$ (that is $\Mo(S,\X)=\{t(\a_1,...,\a_n): t \mbox{ is a term, and } \a_1,...,\a_n\in H\}$).
We usually identify $H$ with $\X$.

\begin{definition}\label {def: zs}
$\zs$, if it exists, is defined to be an E-M set with the following properties: 
\begin{enumerate}
\item If $\X$ is well-ordered, then $\Mo(\zs, \X)$ is well-founded.  \label{def: zs part 1}
\item $\zs$ is {\em cofinal}: For every term $t$, if the following formula is in $\zs$:
\[
\On(t(x_1,...,x_{n-1}))\implies  t(x_1,...,x_{n-1})< x_n,
\]
(where $\On(x)$ is the predicate that says that $x$ is an ordinal).
\item $\zs$ is {\em remarkable}: For every term $t$, if the following formula is in $\zs$:
\begin{multline*}
\On(t(x_1,...,x_n,...,x_{n+m})) \and t(x_1,...,x_n,...,x_{n+m})<x_n \implies \\
  t(x_1,...,x_n,...,x_{n+m}) = t(x_1,...,x_{n-1},x_{n+m+1},...,x_{n+2m+1}).
\end {multline*}
\end{enumerate}
(When we write $\varphi(x_1,...,x_n)$ above, we implicitly assume that $x_1<....<x_n$, and that all the $x_i$'s represent ordinals.)
\end{definition}

If $\zs$ exists, it cannot belong to $L$, and hence its existence is not provable in \ZFC\ (as is not true in \revmathfont{V=L}), but it follows, for example, if we assume that a measurable cardinal exists.

The following two lemmas correspond to \cite[Lemma 2.3]{Dev84} and \cite[Lemma 2.4]{Dev84}.
We include the proofs for completeness.

\begin{lemma}
If $\X$ is a linear ordering without a greatest element,  then $\X$ is cofinal in the ordinals of ${\Mo(\zs,\X)}$.
\end{lemma}
\begin{proof}
Let $x$ be an ordinal in ${\Mo(\zs,\X)}$; we want to find $k\in \X$ with $x<k$.
There is a term $t$ and ordinals $\vec{h}$ in $\X$ such that $x=t(\vec{h})$.
Since $\X$ has no greatest element, there is $k\in \X$, $k>\vec{h}$ (i.e., $\vec{h}=(h_1,...,h_\ell)$ and $h_1<...<h_\ell<k$).
By the cofinality of $\zs$, the formula $t(\vec{h})< k$ is in $\zs$.
So $x<k$.
\end{proof}

Given a linear ordering $\X$ and $\a\in\X$, we let $\X\upto \a$ be the restriction of $\X$ to the set $\{\b\in \X: \b<\a\}$.

\begin{lemma} \label{le: Devlin 2.4}
If $\X$ is a linear ordering, and $\a\in \X$ is a left-limit point (i.e., $\forall \b<_\X\a\ \exists \g\ (\b<_\X\g<_\X\a)$), then
\[
\Mo(\zs,\X\upto \a) \isom (L_\a)^{\Mo(\zs,\X)} .
\]
\end{lemma}
\begin{proof}
Let $\A= \Mo(\zs,\X)$, and $\B$ be the Skolem Hull of $\X\upto \a$ in $\A$.
So $\B \isom \Mo(\zs,\X\upto \a)$.

We start by proving that $\On^\B=\{x\in \On^\A: x<\a\}$.
By the lemma above, we get that $\X\upto \a$ is cofinal in $\On^\B$, and hence $\On^\B\subseteq \{x\in \On^\A: x<\a\}$.
Now, pick $x\in \On^\A$ with $x<\a$.
There is an increasing tuple of ordinals $\vec{k}\vec{h}$ from $\X$, where $\vec{k}<x \leq\vec{h}$, such that $x=t(\vec{k},\vec{h})$.
Since $\a$ is a left-limit and $\vec{k}<\a$, there exists a tuple of ordinals $\vec{l}$ from $\X\upto \a$ with $\vec{k}<\vec{l}<\a$, and with the same number of elements as the tuple $\vec{h}$.
Then, by indiscernibility, $t(\vec{k},\vec{l})\leq \vec{l}$.
And then, by the remarkable property, $x=t(\vec{k},\vec{h})=t(\vec{k},\vec{l})\in \B$.
This finishes the proof that $\On^\B=\{x\in \On^\A: x<\a\}$.

For all $x\in\A$, the least ordinal $\b_x\in \A$ such that $x\in (L_{\b_x})^\A$ is definable in $\A$ by some term $t(x)=\b_x$.
Thus,  for all $x\in \B$, $\b_x\in \B$, and hence $\b_x<\a$.
It follows that $\B\subseteq L_\a^\A$.
Conversely, if $x\in L_\a^\A$, then $x\in L_\b^\A$ for some $\b<\a$.
This implies that $x$ is definable from a finite set of ordinals below $\a$.
(\revmathfont{V=L} implies that if $x\in L_\b$, then $x$ is definable from finitely many ordinals which are $\leq \b$.)
All these ordinals are in $\B$, and hence so is $x$.
We get that $\B=L_\a^\A$.
\end{proof}

\subsection{The fixed-point structure}    \label{ss: fixed point}

Let $\A=\Mo (\zs,\om\cdot\om^*)$, and $H$ be the associated sequence of ordinal indiscernibles in $\A$ or order type $\om\cdot\om^*$.
We claim that
\[
\A \eqSIN  \A'.
\]

Let $\a$ be the least element of the largest copy of $\om$ in $H$.
(That is, we can decompose $H$ as $\om\cdot\om^* + \{\a\} + \om$).
Let $H_0=H\upto \a$, and let $\B= \Mo(\zs, H_0)$, the Skolem hull of $H_0$ in $\A$.
Since $H_0$ is also isomorphic to $\om\cdot\om^*$, by the Ehrenfeucht-Mostowski theorem, we get that $\B\isom\A$.
Since $\zs$ is cofinal and remarkable, and $\a$ is a left-limit point of $H$, Lemmas \ref {le: Devlin 2.4} implies that $\B=L_\a^\A$.

The subset $L_\a^\A\subseteq\A$ is $\Sic_1$-definable in $\A$ using the element $L_\a^\A\in \A$ as a parameter.
Let
\[
K = \{ (e,\bar{b}): e\in \om^\A, \bar{b}\in \B^{<\om}, \B\models\varphi_e(\bar{b})\},
\]
 where $\varphi_e$ is the $e$th $\Sic_1$ formula.
The set $K$ is first-order definable over $\B=L_\a^\A$, and hence belongs to $L_{\a+1}^\A$.
Here we are using that $\A$ is an $\om$-model (i.e., $\om^\A=\om$), so that for all $e$ and $\bar{b}$, 
\[
\B\models \varphi_e(\bar{b}) \iff \A\models \mbox{``}L_\a^\A\models \varphi_e(\bar{b})\mbox{''}.
\]
The reason why $\A$ is an $\om$-model is that there is no term $t(\xbar)$ such that $\zs$ contains the formulas $t(\xbar)\in \om$ and $\underline{n}<t(\xbar)$ for all $n\in \om$.
We know that $\zs$ does not contain all these formulas because,  for instance, $\Mo(\zs,\om)$ is well-founded, and hence an $\om$-model.
(Let us notice that this is all we use from property (\ref{def: zs part 1}) of Definition \ref {def: zs}.)

Using $L_\a^\A$ and $K$ as parameters, we can get an effective interpretation of $\B'$ in $\A$.
Thus $\B'\leqSIN \A$.
Since $\B\isom \A$, we get $\A \eqSIN \A'$.

\

\section{Coding higher-order arithmetic}

\subsection{Background on higher-order arithmetic} \label{ss: HOA}

We let higher-order arithmetic be defined as follows.
The language $L_\om$ has infinitely many sorts, or {\em types}.
There is a first type, that we call $N$, and intends to represent the natural numbers, and then, for every type $\tau$, we have a type $\P(\tau)$ which intends to represent the power set of $\tau$.
We write $\P^n(N)$ for $\P(\P(...(\P(N)...))$ iterated $n$ times.
All variables have a type, and thus, all quantifiers range over the elements of a certain type.
The language has symbols $0,1,+,\times,<$ that apply to elements of the first type $N$, and a binary relation $x\in y$ that can only be used when $x$ is of a certain type $\tau$ and $y$ has type $\P(\tau)$.
The axioms are the ordered-semi-ring axioms for the elements of the first type $N$, plus extensionality for sets of all types, plus the induction axiom for all subsets of $N$ in $\P(N)$, plus  the full comprehension scheme, which for all types $\tau$ and for all formulas $\psi(z^\tau)$ contains the universal closure of the formula  $\exists y^{\P(\tau)}\forall z^\tau(z^\tau\in y^{\P(\tau)} \leftrightarrow \psi(z^\tau))$.
It is not hard to see that a statement can be proved in higher-order arithmetic if and only if it can be proved in full $n$th-order arithmetic for some $n\in \om$:
The reason is that a proof in higher-order arithmetic can only mention finitely many types.

An $\om$-model of the language $L_\om$ is an $L_\om$-structure where the first-order part is standard, i.e.\ $(\om; 0,1,+,\times,<)$.
These $L_\om$-structures are determined by a sequence $(\E^1,\E^2,\E^3,....)$ where $\E^n\subseteq \P^n(\om)$.
Such a structure $\H$ is a model of higher-order arithmetic if and only if for every $n$ and every formula $\varphi(z^n)$ in the language $L_\om$, with parameters from $\H$, and where $z^n$ is a variable of type $\P^n(N)$, we have that
\[
\{F\in \E^n:  \H\models \varphi(F)\}   \in\  \E^{n+1}.
\]

\subsection{Background on generic copies of a structure}

A important tool in our proofs will be the notion of generic copy of a structure  $\A$ introduced by Ash, Knight, Mennasse and Slaman \cite{AKMS89}, and Chisholm \cite{Chi90}.
We refer the reader to \cite[Chapter 10]{AK00} for the basic properties of this forcing notion, that we now quickly review.
(See also \cite[Section 4]{MonND} for an exposition closer to the one we use here.)

We define the forcing notion $\PP$ to be the set of finite one-to-one partial functions from a set of constants $B$ to $A$, the domain of  the structure $\A$.
Diverting from \cite{AK00}, we will set $B=\{b_0,b_1,...\}=\om$, and we will only consider finite functions defined on initial segments of $\om$.
So, we can think of the conditions of $\PP$ as finite tuples of different elements from $\A$ ordered by inclusion.

A generic $G$ for this forcing gives a bijection from $B$ $(=\om)$ to $A$.
By pulling back the relations and functions from $\A$, $G$ defines a structure $\B$ on $B$ which is isomorphic to $\A$.
Let $L$ be the language of $\A$; assume it is a relational language.
Let $\{\phi_0,\phi_1,...\}$ be an enumeration of all the atomic $(L\cup B)$-sentences.
For each $n$, let $\{\phi_0,...,\phi_{k_n-1}\}$ be the subset of atomic $(L\cup B)$-sentences which only use the first $n$ relations from $L$ and the first $n$ constants from $B$.
(Assume these formulas always come first in the listing.)
Let $D(\B)$ be the atomic diagram of $\B$, that is, $D(\B)\in 2^\om$ and $D(\B)(i)=1$ iff $\B\models \phi_i$.
For each $\pbar\in\PP$, we let $D(\pbar)$ be the fragment of $D(\B)$ of lenght $k_ {| \pbar  |} $, determined by $\pbar$.
That is, $D(\pbar)(i)=1$ iff $\A\models \phi_i(\pbar)$, where, in $\phi_i(\pbar)$, each constant $b_j$ is replaced by $\pbar(j)$.
This way, we have that $D(\B)=\bigcup_{\pbar\in G}D(\pbar)$, and that for each $\si\in 2^{k_n}$ there is a quantifier-free formula $\varphi_\si(x_1,...,x_n)$, such that $\A\models\varphi_\si(\pbar)$ if and only if $\si=D(\pbar)$.

Given an infinitary sentence $\varphi$ in the language $L\cup B$, and given $\pbar\in \PP$, the forcing relation $\pbar\forces \varphi$ is defined as usual:
$\pbar\forces \phi_i$ iff $i<k_{|\pbar|}$ and $\A\models \phi_i(\pbar)$;
$\pbar\forces \bigvee_i\psi_i $ iff for some $i$, $\pbar\forces \psi_i$;
and $\pbar\forces \bigwedge_i\psi_i $ iff for each $i$, and each $\qbar\supseteq\pbar$, there exists $\rbar\supseteq\qbar$, such that  $\rbar\forces \psi_i$.
One can then prove that for any sentence $\varphi$, $\B\models \varphi$ if and only if for some $\pbar\in G$, $\pbar\forces \varphi$, provided that $G$ is generic enough.

It is shown in \cite[Lemma 10.6]{AK00} that there is an effective procedure that, given a sentence $\varphi$, returns a formula $\Force_\varphi$ of the same complexity such that $\pbar\forces \varphi\iff \A\models \Force_\varphi(\pbar)$.
We also note that we can uniformly get $\Sic_1$ sentences $\varphi_{e,n}$ in the language $L\cup B$, such that, given a c.e.\ operator $W_e$ and $n\in \om$, $\B\models \varphi_{e,n}$ if and only if $n\in W_e^{D(\B)}$, namely $\varphi_{e,n}\equiv \bigvee_{\si\in 2^{<\om}: n\in W_e^\si} \varphi_\si(b_0,...,b_m)$.
We write this formula $\varphi_{e,n}$ as ``$n\in W_e^\B$.''
(Here $\{W_e:e\in\om\}$ is the standard list of all c.e.\ operators, and by $W_e^\si$, with $\si\in 2^{<\om}$, we mean the set of $n$'s which are enumerated in $W_e^\si$ in less than $|\si|$ steps, using $\si$ as oracle.)
We observe that $\qbar\forces n\in W_e^\B$ if and only if $n\in W_e^{D(\qbar)}$.

The following application of this forcing notion contains many of the ideas that we will use in later proofs.

\begin{lemma}(Knight \cite[Theorem 10.17]{AK00})    \label {lem: knight}
Given $X\subseteq\om$, we have that $X$ is c.e.\ in all copies of $\A$ if and only if there is a tuple $\abar\in A^{<\om}$ and an enumeration-operator $\Theta$ such that 
\[
X=\Theta(\Sic_1\dtp_\A(\abar)).
\]
\end{lemma}
\begin{proof}
The implication from right to left is fairly straightforward.
Any oracle that computes a presentation of $\A$ can enumerate $\Sic_1\dtp_\A(\abar)$ for any $\abar\in A^{<\om}$, and thus can also enumerate $X=\Theta(\Sic_1\dtp_\A(\abar))$.
Let us now consider the other direction.

Suppose that $X$ is c.e.\ in all copies of $\A$.
Let $G$ be an $X$-arithmetically generic filter in $\PP$, and let $\B$ be the associated copy of $\A$.
There is some $e$ such that $X=W_e^\B$.
Thus, there is some $\pbar\in G$ such that $\pbar\forces X=W_e^\B$.
(This means that $\pbar \forces \bigwedge_{n\in X}(n\in W_e^\B) \& \bigwedge_{n\not\in X}\neg(n\in W_e^\B)$.)
So, we claim that for every $n$ 
\[
n\in X \quad\iff\quad \exists \qbar \in\PP\ \left(\qbar \supseteq \pbar\ \ \and\ \  n\in W_e^{D(\qbar)}\right).
\]
The reason is that, if $n\in X$, there is a $\qbar\in G$, $\qbar \supseteq \pbar$ such that $n\in W_e^{D(\qbar)}$; and if $n\not\in X$, then for all $\qbar\supseteq \pbar$ there exists $\rbar\supseteq \qbar$ such that $\rbar\forces \neg(n\in W_e^\B)$ and hence $\qbar\not\forces n\in W_e^\B$.
Let $\Theta$ be the enumeration operator such that, if for some $\si\in 2^{<\om}$ with $n\in W_e^\si$,  the index for the $\Sic_1$ formula
\[
\exists \xbar\ \big(\pbar\xbar\in \PP \and  D(\pbar,\xbar)=\si \big)
\]
belongs to an oracle $Y$, then $n\in \Theta^Y$ (where ``$\qbar\in \PP$'' is the formula that says that all the elements in the tuple $\qbar$ are different).
It follows that $X=\Theta(\Sic_1\dtp_\A(\pbar)).$
\end{proof}

Let $\{\Theta_e:e\in\om\}$ be the standard computable enumeration of all enumeration-operators.
That is, $\Theta_e(Y)$ is the set of all $k\in\om$ for which we have $(u,k)\in W_e$ for some $u$ with $D_u\subseteq Y$, where $D_u$ is the $u$th finite subset of $\om$ in some standard ordering.
For each $e$ and each $\pbar\in A^{<\om}$, let 
\[
T_e(\pbar) = \Theta_{e} (\Sic_1\dtp_\A(\pbar))  \subseteq \om.
\]
Notice that from the uniformity in the  proof above, we get a computable function $f_1\colon\om\to\om$ such that, for each $e$ and each $X\subseteq\om$, 
\[
\pbar\forces X= W_e^\B  \implies   X= T_{f_1(e)}(\pbar) ,
\]
where $\Theta_{f_1(e)}$ is defined by
\[
\Theta_{f_1(e)} (\Sic_1\dtp_\A(\pbar)) = \{m\in\om: (\exists \qbar\supseteq \pbar)\  m\in W_{e}^{D(\qbar)}\} = \{m\in\om: \pbar\not\forces \neg m\in W_e^\B\} .
\]

The next step is to generalize all this to $X\in\P^n(\om)$.

%




%


\subsection{The model $\H_\A$}

\begin{definition}
Given $F\in \P^n(\om)$ and $Y\subseteq \om$, we say that $F$ is {\em uniformly computably enumerable (u.c.e.)} in $Y$ if there exists a $Y$-c.e.\ set $W\subseteq\om^{n}$ such that 
\[
F=\{\{...\{\{\{m\in\om: \la m,i_1,....,i_{n-1}\ra\in W\}: i_1\in\om\}: i_2\in\om\}:...\}: i_{n-1}\in\om\}.
\]
In this case, we say that $W$ {\em codes} $F$.
Given a structure $\A$ and $n\in \om$ we let $\E^n_\A\subseteq \P^n(\om)$ be the set of all $F\in \P^n(\om)$ such that, for every set $X\in\Sp(\A)$, $F$ is u.c.e.\ in $X$.
\end{definition}

So, for instance, $\E^1_\A$ is the set of all subsets of $\om$ which are c.e.\ in all copies of $\A$.
Lemma \ref {lem: knight} says that 
\[
\E^1_\A = \{ T_e(\abar): e\in \om, \abar\in A^{<\om}\}.
\]

\begin{definition}
Given a structure $\A$, we let $\H_\A$ be the $\om$-model determined by
\[
\H_\A=(\E^1_\A, \E^2_\A,\E^3_\A,...),
\]
that is,  $\H_\A$ is the $L_\om$-structure where $(N,0,1,+,\cdot,<)$ is represented by the standard model of arithmetic, and $\P^n(N)$ is interpreted as $\E^n_\A$.
\end{definition}

\begin{theorem}\label{thm: main HA}
If $\Sp(\A)=\Sp(\A')$, then $\H_\A$ is a model of higher-order arithmetic.
\end{theorem}

The proof of this theorem will be completely contained in \Ztwo, and hence Theorem \ref {thm: omega arithmetic} and Corollary \ref {cor: omega arithmetic} follow.
We let the reader verify that the proof goes through in \Ztwo. 
We just notice that $\H_\A$ is a countable structure, and that it can be represented by a countably-coded $\om$-model.

As we mentioned above, all we need to prove is the comprehension axiom scheme.

The next few lemmas are dedicated to characterize $\E^n_\A$ in terms of the types realized in $\A$.
We use $\Sic_m\dtp_\A(\abar)$ to denote the $\Sic_m$-type of $\abar$ in $\A$.

\begin{lemma}
If $F\in \E^2_\A$, then there is a tuple $\abar\in A^{<\om}$ and a uniformly computably infinitary list of $\Sic_3$ formulas $\varphi_\ell(\xbar,\ybar)$, such that 
\[
F=\{ T_\ell(\abar, \bbar): \ell\in\om, \bbar\in A^{<\om},  \A\models \varphi_\ell(\abar, \bbar)\}.
\]
\end{lemma}
\begin{proof}
This proof is somewhat similar to the one of Lemma \ref {lem: knight}.

Since $F\in \E^2_\A$, $F$ is uniformly computably enumerable in the diagram of $\B$ (which, abusing notation, we denote as $\B$).
Then, for some computable sequence $\{e_i:i\in\om\}$, $F$ is coded by $\bigoplus_{i\in\om}W_{e_i}^\B$; that is, $F=\{W_{e_i}^\B: i\in\om\}$.
There exists a condition $\pbar\in \PP$ such that $\pbar\forces F=\{W_{e_i}^\B: i\in\om\}$. 
If 
\[
F=\{X_0,X_1,...\},
\]
this means that $\pbar $ forces that $\forall i\in\om\exists j (W_{e_i}^\B= X_j)$ and $\forall j\exists i (W_{e_i}^\B= X_j)$.

We say that $\qbar$ {\em decides} $W_e^\B$, if for all $k\in\om$, either $\qbar\forces k\not\in W_e^\B$ or $\not\exists \rbar_1\supseteq\qbar\ (\rbar_1\forces  k\not\in W_e^\B)$.
Equivalently, $\qbar$ {decides} $W_e^\B$,  if for all generic extensions $G$ of $\qbar$, $W_e^\B$ is always the same set.
Let $\delta_e(\qbar)$ be the $\Pic_2$ formula that says that $\qbar$ {decides} $W_e^\B$:
\[
\delta_e(\qbar)  \equiv   \bigwedge_{k\in\om} (\forall \rbar\supseteq\qbar (k\not\in W_e^{D(\rbar)})) \vee (\forall \rbar_1\supseteq\qbar\ \exists \rbar_2\supseteq\rbar_1 (k\in W_e^{D(\rbar_2)})).
\]

So, we have that, if for some set $X$, $\qbar\forces W_e^\B=X$, then $\delta_e(\qbar)$ holds.
Conversely, if $\delta_e(\qbar)$ holds, then for some $X$, $\qbar\forces W_e^\B=X$:
Because, for some $\rbar\supseteq\qbar$,  there is an $X$ such that $\rbar\forces W_e^\B=X$, but then $\delta_e(\qbar)$ implies that $\qbar\forces W_e^\B=X$ too.
Recall that this implies that $X= T_{f_1(e)}(\qbar)$.

Let us go back to the fact that $\pbar$ forces that $\forall i\in\om\exists j (W_{e_i}^\B= X_j)$ and $\forall j\exists i (W_{e_i}^\B= X_j)$.
We claim that 
\[
F=\{ T_{f_1(e_i)}(\qbar): \qbar\supseteq\pbar, i\in\om, \delta_{e_i}(\qbar)\}.
\]
Suppose first that $\qbar\supseteq \pbar$ and that $\delta_{e_i}(\qbar)$ holds.
Then for some $j\in \om$ and some $\rbar\supseteq\qbar$, $\rbar\forces W_{e_i}^\B=X_j$.
Therefore, $\qbar\forces W_{e_i}^\B=X_j$ too, and $T_{f_1(e_i)} (\qbar) = X_j$.
So the right-hand side is included in the left-hand side.
Take now $X_j\in F$.
There exists $\qbar\supseteq \pbar$, $\qbar\in G$ and $i\in\om$ such that $\qbar\forces W_{e_i}^\B=X_j$.
Then $\delta_{e_i}(\qbar)$ holds and $T_{f_1(e_i)}(\qbar) = X_j$.
So the left-hand side is included in the right-hand side.

To get the sequence of $\Sic_3$ formulas we wanted, we let  $\varphi_\ell(\abar,\bbar)\equiv \bigvee_{i\in\om: f_1(e_i)=\ell} \delta_{e_i}(\abar,\bbar)$.
\end{proof}

We now want to generalize this proof to $\E^n_\A$.

\begin{definition}
Given $m,n,e\in\om$,  and $\abar\in A^{<\om}$, we define $T_{m,e}^n(\abar) \in \P^n(\om)$ by induction on $n$ as follows.
\begin{eqnarray*}
T_{m,e}^1(\abar)    &=& \Theta_e(\Sic_{m}\dtp_\A(\abar)),  \\
T_{m,e}^{n+1}(\abar) &=&  \{T_{m,j}^n(\abar,\bbar): \bbar\in A^{<\om}, j\in T_{m,e}^1(\abar,\bbar) \}.
\end{eqnarray*}
\end{definition}

Notice that $T^1_{1,e}(\abar)=T_e(\abar)$.

\begin{lemma}
For all $n\geq 1$, 
$ \{T_{1,e}^n(\abar): e\in \om, \abar\in A^{<\om}\} \subseteq\ \E^n_\A$.
\end{lemma}
\begin{proof}
This proof is like the right-to-left proof of Lemma \ref {lem: knight}.
Any oracle that computes a presentation of $\A$ can uniformly enumerate $\Sic_1\dtp_\A(\bbar)$ for any $\bbar\in A^{<\om}$, and can thus enumerate a set coding $T_{1,e}^n(\abar)$.
\end{proof}

Our next goal is to show that, if $\Sp(\A)=\Sp(\A')$, then, for every $n\geq 1$, 
\[
\E^n_\A = \{T_{n^2,e}^n(\abar): e\in \om, \abar\in A^{<\om}\} =  \{T_{m,e}^n(\abar): e,m\in \om, \abar\in A^{<\om}\}.
\]
We have already shown the first equality for the case $n=1$.

Before we prove further results about these families $T^n_{m,e}(\abar)$, let as show how the basic statements about them can be translated to statements about the structure $\A$.
First, given $k, m, e$, note that
\[
k\in T^1_{m,e}(\abar)  \iff \A\models \bigvee_{u: (u,k)\in \Theta_e} \ \left( \bigwedge_{j\in D_u} \varphi_{m,j}(\abar)\right),
\]
where $\{\varphi_{m,j}: j\in \om\}$ is a standard enumeration of all $\Sic_m$ formulas.
We let $\overline{k\in T^1_{m,e}(\xbar)}$ denote this $\Sic_{m}$ formula $\bigvee_{u: (u,k)\in \Theta_e} \  \bigwedge_{j\in D_u} \varphi_{m,j}(\xbar)$.
Note that this is now a formula in the language of $\A$.

We now define the following $\Pic_{m+1}$ formula, where $m=\max\{m_1,m_2\}$.
\[
\overline {T_{m_1,e_1}^{1}(\xbar_1)= T_{m_2,e_2}^{1}(\xbar_2)}  \ \ \equiv \ \ \bigwedge_{k\in\om}\left(  \overline{k\in T^1_{m_1,e_1}(\xbar_1)} \leftrightarrow \overline{k\in T^1_{m_2,e_2}(\xbar_2)}\right).
\]

We now want to define a formula $\overline{T_{m_1,e_1}^{n}(\xbar_1)\subseteq T_{m_2,e_2}^{n}(\xbar_2)}$ in the language of $\A$ such that  for all $\abar_1,\abar_2\in A^{<\om}$,
\[
T_{m_1,e_1}^{n}(\abar_1) \subseteq  T_{m_2,e_2}^{n}(\abar_2)  \iff \A\models \overline{T_{m_1,e_1}^{n}(\abar_1)\subseteq T_{m_2,e_2}^{n}(\abar_2)}.
\]
We let

\begin{multline*}
\overline{T_{m_1,e_1}^{n}(\xbar_1) \subseteq  T_{m_2,e_2}^{n}(\xbar_2)}\  \ \equiv\ \  \bigwedge_{i\in\om}\forall \ybar_1 \bigg( \overline{i\in T^1_{m_1,e_1}(\xbar_1,\ybar_1)} \implies       \\ 
\bigvee_{j\in\om} \exists \ybar_2 \left( \overline{j\in T^1_{m_2,e_2} (\xbar_2,\ybar_2)} \ \  \and \ \ \overline{T_{m_1,i}^{n-1}(\xbar_1,\ybar_1)= T_{m_2,j}^{n-1}(\xbar_2,\ybar_2)} \right)\bigg).
 \end{multline*}
It is not too hard to see, using induction on $n$, that $\overline{T_{m_1,e_1}^{n}(\xbar_1)\subseteq T_{m_2,e_2}^{n}(\xbar_2)}$ is a $\Pic_{m+2n-1}$  formula where $m=\max\{m_1,m_2\}$.

Of course, we define the $\Pic_{m+2n-1}$ formula $\overline{T_{m_1,e_1}^{n}(\xbar_1) =  T_{m_2,e_2}^{n}(\xbar_2)}$ using both inclusions in the obvious way.

\begin{lemma}
For every $n\geq 1$, $\E^n_\A \subseteq \{T_{n^2,e}^n(\abar): e\in \om, \abar\in A^{<\om}\}$.
\end{lemma}
\begin{proof}

By induction on $n$, we prove that there is a computable function $f_n$ such that, for all $F\in \P^n(\om)$, $e\in\om$, $\pbar, \pbar_1,\pbar_2\in \PP$
\[
\pbar \forces `` W_e^{\B} \mbox{ codes } F \mbox{''}  
	\quad\implies\quad
T^n_{n^2, f_n(e)}(\pbar) = F,
\]
and
\[
\pbar_1 \subseteq \pbar_2 
        \quad\implies\quad
T^n_{n^2, f_n(e)}(\pbar_1) \supseteq T^n_{n^2, f_n(e)}(\pbar_2).
\]
Notice that we have already proved the case $n=1$ in Lemma \ref {lem: knight}.

Suppose we have already defined such a function $f_n$, and we now want to define $f_{n+1}$.
First, for each $e$, let $\delta_e$ be the following $\Pic_{n^2+2n-1}$ formula:
\[
\delta_e(\pbar)  \ \   \equiv\ \   \forall \qbar\supseteq\pbar\  \big(\overline{T^n_{n^2, f_n(e)}(\pbar) = T^n_{n^2, f_n(e)}(\qbar)}\big).
\]
Notice that if for some $F\in \P^n(\om)$, $\pbar \forces \left(W_e^{D(\A)} \mbox{ codes } F\right) $, then $\delta_e(\pbar)$ holds by the induction hypothesis.
Also note that $n^2+2n-1< (n+1)^2$, so $\delta_e$ is $\Sic_{(n+1)^2}$.

Now, suppose that $F=\{F_1,F_2,...\} \in \P^{n+1}(\om)$ and that
\[
\pbar \forces W_e^{\B} \mbox{ codes } F,
\]
where $W_e^{\B}=\bigoplus_iW_{e_i}^\B$.
That means that
\[
\pbar \forces \left(\forall i\exists j (W_{e_i}^\B \mbox{ codes } F_j)\right)  \& \left(\forall j\exists i (W_{e_i}^\B \mbox{ codes } F_j)\right).
\]
This implies that for each $i$ and each $\qbar\supseteq\pbar$, there exists $j$ and $\qbar_1\supseteq \qbar$ such that $\qbar_1\forces W_{e_i}^\B \mbox{ codes } F_j$.
Now, by the inductive hypothesis, we have that $T^n_{n^2, f_n(e_i)}(\qbar_1) = F_j$.
If we had $\delta_{e_i}(\qbar)$ we would have that $T^n_{n^2, f_n(e_i)}(\qbar) = F_j$ too.
So, for each $i$ and each $\qbar\supseteq\pbar$ with $\delta_{e_i}(\qbar)$, we have that for some $j$, $T^n_{n^2, f_n(e_i)}(\qbar) = F_j$.
On the other hand, for each $j$, there exists $i$ and  $\qbar\supseteq\pbar$ such that $\qbar\forces W_{e_i}^\B \mbox{ codes } F_j$.
For this $\qbar$, $\delta_{e_i}(\qbar)$ holds.
It follows that
\[
F   = \{ T^n_{n^2, f_n(e_i)}(\qbar): i\in \om, \qbar\supseteq \pbar, \delta_{e_i}(\qbar)\}.
\]
Now, let $\Theta_{f_{n+1}(e)}$ be an enumeration-operator such that $k\in \Theta_{f_{n+1}(e)}(\Sic_{(n+1)^2}\dtp_\A(\qbar))$ if and only if for some $i$, $f_n(e_i)=k$ and $\A\models \delta_{e_i}(\qbar)$. 
Then $F=T^{n+1}_{(n+1)^2, f_{n+1}(e)}(\pbar)$ as wanted.
It is easy to see from the definition that $\pbar_1 \subseteq \pbar_2 \implies T^{n+1}_{(n+1)^2, f_{n+1}(e)}(\pbar_1) \supseteq T^{n+1}_{(n+1)^2, f_{n+1}(e)}(\pbar_2)$.
\end{proof}

\begin{corollary}
If $\Sp(\A)=\Sp(\A')$ then, for every $n$, 
\[
\E^n_\A = \{T_{n^2,e}^n(\abar): e\in \om, \abar\in A^{<\om}\} =  \{T_{m,e}^n(\abar): e,m\in \om, \abar\in A^{<\om}\}.
\]
\end{corollary}
\begin{proof}
The inclusion of the first set in the second set follows from the lemma above.
The inclusion of the second set in the third set is obvious.

The hypothesis $\Sp(\A)=\Sp(\A')$ is only used to prove the inclusion $\{T_{m,e}^n(\abar): e,m\in \om, \abar\in A^{<\om}\}\subseteq \E^n_\A$.
If $Y\in\Sp(\A)$, then $Y\in \Sp(\A^{(m-1)})$ and hence $Y$ computes a copy of $\A^{(m-1)}$ and can uniformly enumerate $\Sic_{m}\dtp_\A(\bbar)$ for all $\bbar$ in this copy.
Thus  $Y$ can enumerate $T_{m,e}^n(\abar)$.
\end{proof}

We are now ready to prove that $\H_\A$ satisfies full comprehension. 
Add to $L_\om$, the language  of higher-order arithmetic, a new constant symbol $T^n_{m,e}(\xbar)$ of type $\P^n(N)$ for each $n,m,e\in\om$ and tuple of variables $\xbar$.
(We use $x$ and $y$ for variables that range over elements of $\A$, and $z^n$ for variables in $L_\om$ of type $\P^n(N)$.)
We now define a computable transformation that takes a sentence $\psi(T_{m_1,e_1}^{n_1}(\xbar_1), ... ,T_{m_k,e_k}^{n_k}(\xbar_k))$ in this language, and returns a $\Pi^c_{<\om}$ formula $\overline{\psi_{\bar{m},\bar{e}} (\xbar_1,...,\xbar_k)}$ in the language of $\A$ such that, for all $\abar_1,...,\abar_k\in A^{<\om}$,
\[
\H_\A\models   \psi (T_{m_1,e_1}^{n_1}(\abar_1), ... ,T_{m_k,e_k}^{n_k}(\abar_k))   
				\quad \iff\quad
\A\models \overline{\psi_{\bar{m},\bar{e}}(\abar_1,....,\abar_k)},
\]
where $\bar{m}=(m_1,...,m_k)$ and $\bar{e}=(e_1,...,e_k)$.

We have already defined
\begin{itemize}
\item $\overline{k\in T^1_{m,e}(\xbar)}$,
\item $\overline{T_{m_1,e_1}^{n}(\xbar_1) \subseteq  T_{m_2,e_2}^{n}(\xbar_2)}$, and
\item $\overline{T_{m_1,e_1}^{n}(\xbar_1) =  T_{m_2,e_2}^{n}(\xbar_2)}$,
\end{itemize}
with the desired property.
We now add to the list

\begin{itemize}
\item   $\overline{T_{m_1,e_1}^{n}(\xbar_1)\in T_{m_2,e_2}^{n+1}(\xbar_2)}\ \equiv \    \bigvee_j \exists \ybar \left( \overline{j\in T^1_{m_2,e_2} (\xbar_2,\ybar)} \ \  \and \ \ \overline{T_{m_1,e_1}^{n}(\xbar_1)= T_{m_2,j}^{n}(\xbar_2,\ybar)}\right)$
\item $\overline{ \varphi\vee\psi} = \overline{\varphi}\vee \overline{\psi}$.
\item $\overline{ \neg\psi} = \neg \overline{\psi}$.
\item $\overline{ \exists z^n \psi(z^n)} =\bigvee_{e\in\om} \exists \xbar\ \overline{\psi(T^n_{n^2,e}(\xbar))}$.
\end{itemize}

This transformation allows us to prove comprehension in $\H_\A$ and finish the proof of Theorem \ref {thm: main HA}.

Suppose we have a sentence $\psi(\bar{T},z^n)$ in the language of $\H_\A$ with parameters $\bar{T}=$ $T_{m_1,e_1}^{n_1}(\abar_1)$, ... ,$T_{m_k,e_k}^{n_k}(\abar_k)\in\H_\A$, and a free variable $z^n$ of type $\P^n(N)$.
We need to show that the set
\[
\{ G\in \E^n_\A: \H_\A\models \psi(\bar{T},G)\}
\]
belongs to $\E^{n+1}_\A$.
(Let $\abar=(\abar_1,...,\abar_k)$.)
Then, we have that,
\[ 
\{G\in  \E^n_\A: \H_\A\models \psi(\bar{T},G)\} =  \{T_{n^2,e}^n(\abar,\bbar): e\in\om, \bbar\in A^{<\om},\   \A\models \overline{\psi_{\bar{m},\bar{e}}(\abar,\bbar)} \},
\]  
where $\bar{m}=(m_1,...,m_k,n^2)$ and $\bar{e}=(e_1,...,e_k,e)$.
Let $m$ be such that $\overline{\psi_{\bar{m},\bar{e}}}$ is $\Sic_m$.
It follows that for the appropriate index $j$, $\{G\in  \E^n_\A: \H_\A\models \psi(\bar{T},G)\}= T_{m,j}^{n+1}\in \E^{n+1}_\A$.

This concludes the proof that $\H_\A$ satisfies the full comprehension axiom, and hence that it is a model of higher-order arithmetic.


\begin{thebibliography}{AKMS89}

\bibitem[AK00]{AK00}
C.J. Ash and J.~Knight.
\newblock {\em Computable Structures and the Hyperarithmetical Hierarchy}.
\newblock Elsevier Science, 2000.

\bibitem[AKMS89]{AKMS89}
Chris Ash, Julia Knight, Mark Manasse, and Theodore Slaman.
\newblock Generic copies of countable structures.
\newblock {\em Ann. Pure Appl. Logic}, 42(3):195--205, 1989.

\bibitem[Bal06]{Bal06}
V.~Baleva.
\newblock The jump operation for structure degrees.
\newblock {\em Arch. Math. Logic}, 45(3):249--265, 2006.

\bibitem[Chi90]{Chi90}
John Chisholm.
\newblock Effective model theory vs.\ recursive model theory.
\newblock {\em J. Symbolic Logic}, 55(3):1168--1191, 1990.

\bibitem[Dev84]{Dev84}
Keith~J. Devlin.
\newblock {\em Constructibility}.
\newblock Perspectives in Mathematical Logic. Springer-Verlag, Berlin, 1984.

\bibitem[DJ94]{DJ94}
Rod Downey and Carl~G. Jockusch.
\newblock Every low {B}oolean algebra is isomorphic to a recursive one.
\newblock {\em Proc. Amer. Math. Soc.}, 122(3):871--880, 1994.

\bibitem[EP70]{EP70}
H.~B. Enderton and Hilary Putnam.
\newblock A note on the hyperarithmetical hierarchy.
\newblock {\em J. Symbolic Logic}, 35:429--430, 1970.

\bibitem[Har68]{Har68}
J.~Harrison.
\newblock Recursive pseudo-well-orderings.
\newblock {\em Transactions of the American Mathematical Society},
  131:526--543, 1968.

\bibitem[Kal09]{Kal09}
I.~Sh. Kalimullin.
\newblock Relations between algebraic reducibilities of algebraic systems.
\newblock {\em Izv. Vyssh. Uchebn. Zaved. Mat.}, 53(6):71--72, 2009.

\bibitem[Khi04]{Khi04}
A.~N. Khisamiev.
\newblock On the {E}rshov upper semilattice {$L_E$}.
\newblock {\em Sibirsk. Mat. Zh.}, 45(1):211--228, 2004.

\bibitem[Mon09]{MonJumpStr}
Antonio Montalb\'an.
\newblock Notes on the jump of a structure.
\newblock {\em Mathematical Theory and Computational Practice}, pages 372--378,
  2009.

\bibitem[Mon10]{MonND}
Antonio Montalb\'an.
\newblock Coding and definability in computable structures.
\newblock Notes from a course at Notre Dame University to be published in the
  NDJFL, 2010.

\bibitem[Mor04]{Mor04}
A.~S. Morozov.
\newblock On the relation of {$\Sigma$}-reducibility between admissible sets.
\newblock {\em Sibirsk. Mat. Zh.}, 45(3):634--652, 2004.

\bibitem[Puz09]{Puz09}
V.~G. Puzarenko.
\newblock On a certain reducibility on admissible sets.
\newblock {\em Sibirsk. Mat. Zh.}, 50(2):415--429, 2009.

\bibitem[Ric77]{Ric77}
Linda Richter.
\newblock {\em Degrees of unsolvability of models}.
\newblock PhD thesis, University of Illinois at Urbana-Champaign, 1977.

\bibitem[Sos07]{Sos07}
Alexandra Soskova.
\newblock A jump inversion theorem for the degree spectra.
\newblock In {\em Proceeding of CiE 2007}, volume 4497 of {\em Lecture Notes in
  Comp. Sci.}, pages 716--726. Springer-Verlag, 2007.

\bibitem[SS09]{SS09}
Alexandra~A. Soskova and Ivan~N. Soskov.
\newblock A jump inversion theorem for the degree spectra.
\newblock {\em J. Logic Comput.}, 19(1):199--215, 2009.

\bibitem[Stu]{StuEMU}
A.~I. Stukachev.
\newblock Effective model theory via the $\sigma$-definability approach.
\newblock To appear in the proccedings of EMU.

\bibitem[Stu07]{Stu07}
A.~I. Stukachev.
\newblock Degrees of presentability of models. {I}.
\newblock {\em Algebra Logika}, 46(6):763--788, 793--794, 2007.

\bibitem[Stu08]{Stu08}
A.~I. Stukach{\"e}v.
\newblock On degrees of presentability of models. {II}.
\newblock {\em Algebra Logika}, 47(1):108--126, 131, 2008.

\bibitem[Stu10]{Stu10}
A.~I. Stukachev.
\newblock A jump inversion theorem for the semilattices of {S}igma-degrees
  [translation of mr2586684].
\newblock {\em Siberian Adv. Math.}, 20(1):68--74, 2010.

\end{thebibliography}

\end{document}